\documentclass{amsart}
\usepackage{amsmath}
  \usepackage{paralist}
  \usepackage{graphics} 
  \usepackage{epsfig} 
\usepackage{graphicx}  \usepackage{epstopdf}
 \usepackage[colorlinks=true]{hyperref}

\hypersetup{urlcolor=blue, citecolor=red}

\newtheorem*{main}{Theorem}
\newtheorem{lemma}{Lemma}

\theoremstyle{definition}

\newtheorem{remark}{Remark}

\newcommand{\ep}{\varepsilon}

\newcommand{\bs}{\boldsymbol}

\newtheorem*{defin}{Definition}

\title[] 
      {On the existence of weak solutions of a thermistor system with $p$-Laplacian type equation: the unsteady case}

\author[]
{}


 \email{jnaumann@math.hu-berlin.de}

\begin{document}
\maketitle

\centerline{\scshape Joachim Naumann}
\medskip
{\footnotesize
 \centerline{Department of Mathematics, Humboldt University Berlin}
   \centerline{Unter den Linden 6, 10099 Berlin, Germany}
   } 

\bigskip


\begin{abstract}
Let $\Omega\subset\mathbb{R}^n$ ($n=2$ or $n=3$) be a bounded domain. We consider the thermistor system
\[
\text{(1)}\quad \nabla\cdot \boldsymbol{J}=0,\qquad \text{(2)}\quad \frac{\partial u}{\partial t}+\nabla\cdot\boldsymbol{q}=f(x,t,u,\nabla\varphi)\;\text{ in }\; \Omega\times\,]\,0,T\,[\,, 
\]
where (1) is a $p$-Laplace type equation for $\varphi$ ($u=$~temperature, $\varphi=$~electrostatic potential). We prove the existence of a weak solution $(\varphi,u)$ of (1)--(2) under mixed boundary conditions for $\varphi$, and a Robin boundary condition and an initial condition for $u$.
\end{abstract}

\section{Introduction}

Let $\Omega\subset\mathbb{R}^n$ ($n=2$ or $n=3$) be a bounded domain with Lipschitz boundary $\partial\Omega$, and set $Q_T=\Omega\times\,]\,0,T\,[\,$ ($0<T<+\infty$).\par
Let $\bs{J}$ and $\bs{q}$ denote the electric current field density and the heat flux, respectively, of a thermistor occupying the domain $\Omega$ under unsteady operating conditions. Then the balance equations for the electric current and the heat flow within the thermistor material are the following two PDEs
\[
 \nabla\cdot\boldsymbol{J}=0,\quad \frac{\partial u}{\partial  t}+\nabla\cdot\bs{q}=f(x,t,u,\nabla\varphi)\quad\text{in }\; Q_T,
\]
where $\varphi=\varphi(x,t)$ and $u=u(x,t)$ represent the electrostatic potential and the temperature, respectively 
(see, e.g., \cite[Chap. 8]{25}).\par
We make the following constitutive assumptions on $\bs{J}$ and $\bs{q}$
\[
 \bs{J}=\sigma\big(u,|\bs{E}|\big)\bs{E}\quad\text{ Kirchhoff's law},\quad  \bs{q}=-\kappa(u)\nabla u\quad\text{ Fourier's law},
\]
where
\begin{align*}
 \bs{E}&=-\nabla\varphi\quad\text{ density of the electric field},\\
 \sigma&=\sigma\big(u,|\bs{E}|\big)\quad\text{ electrical conductivity},\\
 \kappa&=\kappa(u)\quad\text{ thermal conductivity}.
\end{align*}
With these notions the above system of PDEs takes the form
\begin{align}\label{1}
 -\nabla\cdot\big(\sigma\big( u,|\nabla\varphi|\big)\nabla\varphi\big)&=0\quad\text{ in }\;Q_T,\\
 \label{2}\frac{\partial u}{\partial t}-\nabla\cdot\big(\kappa(u)\nabla u\big)&=f(x,t,u,\nabla\varphi)\quad\text{ in }\;Q_T.
\end{align}
The function $f=f(x,t,u,\nabla\varphi)$ represents a heat source that will include the Joule heat $\bs{J}\cdot\bs{E}$ as special case. 
\par
We supplement system (\ref{1})--(\ref{2}) by boundary conditions for $\varphi$ and $u$, and an initial condition for $u$. Without any further reference, throughout the paper we assume\vspace*{1mm}
\[
 \partial\Omega=\Gamma_D\cup\Gamma_N\;\text{ disjoint}, \quad\Gamma_D\;\text{ non-empty, open}.
\]
Define
\[
 \Sigma_D=\Gamma_D\times\,]0,T\,[\,,\quad\Sigma_N=\Gamma_N\times\,]\,0,T\,[\,.
\]
We then consider the conditions
\begin{equation}\label{3}
\varphi=\varphi_D \ \text{ on } \ \Sigma_D,\quad \bs{J}\cdot\bs{n}=0 \ \text{ on } \ \Sigma_N,
\end{equation}
\begin{equation}
\label{4}
\bs{q}\cdot\bs{n}= g(u-h) \ \text{ on } \ \partial\Omega\times\,]\,0,T\,[\,,
\end{equation}
\begin{equation}
\label{5}
u=u_0 \ \text{ in } \ \Omega\times\{0\}
\end{equation}\vspace*{-3mm}

\noindent
($\bs{n}=$~unit outward normal to $\partial\Omega$). The first condition in (\ref{3}) means that there is an applied voltage $\varphi_D$ along $\Sigma_D$, whereas the second condition characterizes electrical insulation of the thermistor along $\Sigma_N$. The Robin boundary condition (\ref{4})\footnote{This boundary condition is also called Newton's cooling law.} means that the flux of heat through $\partial\Omega\times\,]\,0,T\,[\,$ is proportional to the temperature difference $u-h$, where $g$ denotes the thermal conductivity of the surface $\partial\Omega$ of the thermistor, and $h$ represents the ambient temperature (cf. \cite{8}, \cite{10}, \cite{15}, \cite{24} (nonlinear boundary conditions)).\hfill$\square$
\par
We consider the following prototype for electrical conductivities $\sigma$ in (\ref{1}). Let $\sigma_0:\mathbb{R}\to\mathbb{R}_+$\footnote{$\mathbb{R}_+=[0,+\infty\,[\,$.} be a continuous function such that 
\[
 0<\sigma_*\le\sigma_0(u)\le\sigma^*\quad\forall\; u\in\mathbb{R}\quad(\sigma_*,\sigma^*=\mathrm{const}).
\]
Let $\delta=\mathrm{const}>0$ and let $1<p<+\infty$. We consider
\begin{equation}\label{6}
 \sigma=\sigma\big(u,|\xi|\big)=\sigma_0(u)\big(\delta+|\xi|^2\big)^{(p-2)/2},\quad\xi\in\mathbb{R}^n.
\end{equation}
Here, the factor $\sigma_0(u)$ describes the thermal dependence of the electrical conductivity $\sigma$ of the thermistor material. We obtain
\[
 \bs{J}=\sigma\big(u,|\bs{E}|\big)\bs{E}=-\sigma_0(u)\big(\delta+|\nabla\varphi|^2\big)^{(p-2)/2}\nabla\varphi
\]
and equ. (\ref{1}) is of $p$-Laplace type 
\[
 -\nabla\cdot\big(\sigma_0(u)\big(\delta+|\nabla\varphi|^2\big)^{(p-2)/2}\nabla\varphi\big)=0.
\]

If $p=2$ and $f=\bs{J}\cdot\bs{E}$ (Joule heat), then (\ref{1})--(\ref{2}) represent the well-known thermistor system (see, e.g., \cite{1}, \cite{9}).\par
To make things clearer, let $I=|\bs{J}|$ and $V=|\bs{E}|$ denote current and voltage, respectively, in an electrical conductor. With $\sigma$ as in (\ref{6}) we obtain the current-voltage characteristic
\begin{equation}
 \label{7} 
 I=I(u,V)=\sigma_0(u)(\delta+V^2)^{(p-2)/2}V.
\end{equation}
If $1<p\le 2$ ($p$ ``near to 1''), then this characteristic  describes approximately the current-voltage relations of transistors (cf., e.g., \cite{16}, \cite[Chap. 6.2.2]{27}). In particular, if $p=1$, then (\ref{7}) is widely used to model the effect of saturation of current under high electric fields in certain transistors. For details see, e.g., \cite[Chap. 2.5]{21}.

\begin{remark}\label{r1}
(The case $2\le p<+\infty$). In \cite{13} (formula (\ref{1}), $\alpha\ge 1$; $p=\alpha+1$ in our notation), the authors consider current-voltage characteristics of the form 
\begin{equation}\label{8}
 I=I(u,V)=\sigma_0(u)V^{p-2},\quad 2\le p<+\infty
\end{equation}
for modeling thermistor-like self-heating effects in organic semiconductors. These characteristics can be approximated by (\ref{7}) for sufficiently small $\delta>0$.\par
For the steady case of (\ref{1})--(\ref{4}) and coefficients $\sigma_0=\sigma_0(x,u)$ and $2<p<+\infty$ in (\ref{8}), the existence of weak solutions for the case of two dimensions has been proved for the first time in \cite{18}. This result was extended to the case of measurable exponents $2\le p(x)<+\infty$ ($x\in\Omega$) in \cite{14}. An extension of the latter result has been recently presented in \cite{7}.\hfill$\square$
\end{remark}

We present a prototype for functions $f=f(x,t,u,\xi)$ on the right hand side of (\ref{2}). For $(x,t,u,\xi)\in Q_T\times\mathbb{R}\times\mathbb{R}^n$ we consider
\begin{equation}\label{9}
 f(x,t,u,\xi)=\eta\big(x,t,u,-a(u,-\xi)\xi\big)\sigma\big(u,|\xi|\big)|\xi|^2,
\end{equation}
where $\sigma=\sigma\big(u,|\xi|\big)$ is as in (\ref{6}) and
\[
 \left\{\begin{array}{l}
         \eta=\eta(x,t,u,\hat{\xi}):Q_T\times\mathbb{R}\times\mathbb{R}^n\to\mathbb{R}_+ \ \text{ is Carath\'eodory},\\[1mm]
         0\le\eta(x,t,u,\hat{\xi})\le\eta_1=\mathrm{const} \quad \forall\;(x,t,u,\hat{\xi})\in Q_T\times\mathbb{R}\times\mathbb{R}^n,\\[1mm]
         a=a(u,\xi):\mathbb{R}\times\mathbb{R}^n\to\mathbb{R} \ \text{ is continuous}.
        \end{array}\right.
\]
Writing $\xi=-\bs{E}$ \footnote{Recall $\bs{E}=-\nabla\varphi$; $\bs{J}=\sigma\big(u,|\bs{E}|\big)\bs{E}$.} we obtain
\[
 f(x,t,u,-\bs{E})=\eta\big(x,t,u,a(u,\bs{E})\bs{E}\big)\bs{J}\cdot\bs{E}.
\]
The condition for $\eta$ and $a$ can be specified in several ways, e.g., $a(u,\xi)=\sigma\big(u,|\xi|\big)$. Then $\eta$ may be considered as depending on $\bs{J}$. In particular, if $0<\eta(x,t,u,\bs{J})<1$, then the source term $f(x,t,u,-\bs{E})$ in (\ref{2}) models a loss of Joule heat (cf. \cite{18} for more details).\hfill $\square$
\bigskip

Our paper is organized as follows
\begin{itemize}
 \item[2.] Weak formulation of (\ref{1})--(\ref{5}). Statement of the main result
 \item[3.] Proof of the main result
 \begin{itemize}
  \item[3.1] Existence of an approximate solution
  \item[3.2] A-priori estimates
  \item[3.3] Passage to the limit $\ep\to 0$
 \end{itemize}
 \end{itemize}\par
 References
 

\section{Weak formulation of (\ref{1})--(\ref{5}). Statement of the main result}

We introduce the notations which will be used in what follows.\par
By $W^{1,p}(\Omega)$ ($1\le p<+\infty$) we denote the usual Sobolev space. Define
\[
 W_{\Gamma_D}^{1,p}(\Omega)=\big\{ v\in W^{1,p}(\Omega); v=0 \ \text{ a.e. on } \ \Gamma_D\big\}.
\]
This space is a closed subspace of $W^{1,p}(\Omega)$. Throughout the paper, we consider $W_{\Gamma_D}^{1,p}(\Omega)$ equipped with the norm
\[                                                                                                                                             |v|_{W^{1,p}}=\left(\,\int\limits_\Omega|\nabla v|^p dx\right)^{1/p}.                                                                                                                                            \]
Let $X$ denote a real normed space with norm $|\cdot|_X$ and let $X^*$ be its dual space. By $\langle x^*,x\rangle_X$ we denote the dual pairing between $x^*\in X^*$ and $x\in X$. The symbol $L^p(0,T,X)$ ($1\le p\le+\infty$) stands for the vector space of all strongly measurable mappings $u:\,]\,0,T\,[\,\to X$ such that the function $t\mapsto\big|u(t)\big|_X$ is in $L^p(0,T)$ (cf. \cite[Chap. III, \S3; Chap. IV, \S3]{4}, \cite[App.]{5}, \cite[Chap. 1]{11}). For $1\le p<+\infty$, the spaces $L^p\big(0,T;L^p(\Omega)\big)$ and $L^p(Q_T)$ are linearly isometric. Therefore, in what follows we identify these spaces. \par
Let $H$ be a real Hilbert space with scalar product $(\cdot,\cdot)_H$ such that $X\subset H$ densely and continuously. Identifying $H$ with its dual space $H^*$ via Riesz' Representation Theorem, we obtain the continuous embedding $H\subset X^*$ and 
\[
 \langle h,x\rangle_X=(h,x)_H\quad\forall\; h\in H, \ \forall\; x\in X.
\]
Given any $u\in L^1(0,T;X)$ we identify this function with a function in $L^1(0,T;X^*)$ and denote it again by $u$. If there exists $U\in L^1(0,T;X^*)$ such that 
\[
 \int\limits_0^T u(t)\alpha'(t)dt\mathop{=}\limits^{\mathrm{in } X^*}-\int\limits_0^T U(t)\alpha(t)dt\quad\forall\;\alpha\in C_c^\infty(\,]\,0,T\,[\,),
\]
then $U$ will be called derivative of $u$ in the sense of distributions from $\,]\,0,T\,[\,$ into $X^*$ and denoted by $u'$ (see \cite[App.]{5}, \cite[Chap. 21]{11}).\hfill$\square$
\par
Let $1<p<+\infty$ be fixed. We make the following assumptions on the coefficients $\sigma$, $\kappa$ and the right hand side $f$ in (\ref{1})--(\ref{2}):
\[
\begin{array}{l}
\text{(H1)}\qquad\left\{\begin{array}{l}
           \sigma:\mathbb{R}\times\mathbb{R}_+\to\mathbb{R}_+ \ \text{ is continuous},\\[1mm]
           c_1\tau^p-c_2\le\sigma(u,\tau)\tau^2,\; 0\le\sigma(u,\tau)\le c_3(1+\tau^2)^{(p-2)/2}\\[1mm]
           \forall\;(u,\tau)\in\mathbb{R}\times\mathbb{R}_+, \text{ where } c_1,c_3=\mathrm{const}>0\text{ and } c_2=\mathrm{const}\ge0;
                        \end{array}\right.\\[7mm]
 \text{(H2)}\qquad \left\{\begin{array}{l}
                           \kappa:\mathbb{R}\to\mathbb{R}_+ \ \text{ is continuous},\\[1mm]
                           0 <\kappa_0\le\kappa(u)\le\kappa_1 \;\: u\in\mathbb{R}, \text{ where } \ \kappa_0,\kappa_1=\mathrm{const},
                          \end{array}\right.
                          \end{array}
\]
and the natural growth condition (with respect to (H1))
\[
 \text{(H3)}\qquad\left\{\begin{array}{l}
                   f:Q_T\times\mathbb{R}\times\mathbb{R}^n\to\mathbb{R}_+ \ \text{ is Carath\'eodory}, \\[1mm] 
                   0\le f(x,t,u,\xi)\le c_4\big(1+|\xi|^p\big)\\[1mm]
                   \forall\; (x,t,u,\xi)\in Q_T\times\mathbb{R}\times\mathbb{R}^n,\text{ where } \ c_4=\mathrm{const} >0.\qquad\qquad\quad
                        \end{array}\right.
\]
It is readily seen that (H1) and (H3) are satisfied by the prototypes for $\sigma$ and $f$ we have considered in Section 1.

\begin{defin}
 Assume (H1)--(H3) and suppose that the data in (\ref{3})--(\ref{5}) satisfy
 \begin{equation}\label{10}
  \varphi_D\in L^p\big(0,T;W^{1,p}(\Omega)\big);
 \end{equation}
\begin{equation}
 \label{11}
 g=\mathrm{const}, \quad h=\mathrm{const};
\end{equation}
\begin{equation}
 \label{12}
 u_0\in L^1(\Omega).
\end{equation}
The pair
\[
(\varphi,u)\in L^p\big(0,T;W^{1,p}(\Omega)\big)\times L^q\big(0,T;W^{1,q}(\Omega)\big)\quad \Big(1<q<\frac{n+2}{n+1}\Big)
\]
is called weak solution of (\ref{1})--(\ref{5}) if
\begin{equation}\label{13}
 \exists\; u'\in L^1\big(0,T;\big(W^{1,q'}(\Omega)\big)^*\big);
\end{equation}
\begin{equation}\label{14}
 \int\limits_{Q_T}\sigma\big(u,|\nabla\varphi|\big)\nabla\varphi\cdot\nabla\zeta dxdt=0\quad\forall\;\zeta\in L^p\big(0,T;W_{\Gamma_D}^{1,p}(\Omega)\big);
\end{equation}
\begin{equation}\label{15}
\left.\begin{array}{l}
        {\displaystyle\int\limits_0^T\langle u',v\rangle_{W^{1,q'}}dt+\int\limits_{Q_T}\kappa(u)\nabla u\cdot\nabla v\, dxdt+g\int\limits_0^T\int\limits_{\partial\Omega}(u-h)v\,d_x Sdt}\\
        ={\displaystyle \int\limits_{Q_T}f(x,t,u,\nabla\varphi)v\,dxdt\quad\forall\; v\in L^\infty\big(0,T;W^{1,q'}(\Omega)\big);}
       \end{array}\right\}
\end{equation}
\begin{equation}\label{16}
 \varphi=\varphi_D\quad\text{ a.e. on } \ \Sigma_D;
\end{equation}
\begin{equation}\label{17}
 u(0)=u_0\quad\text{ in } \ \big(W^{1,q'}(\Omega)\big)^*.
\end{equation}
\end{defin}
\medskip

\noindent
The condition $1<q<\frac{n+2}{n+1}$ is standard for weak solutions of parabolic equations with right hand side in $L^1$.\par
We notice that the function $f(\cdot,\cdot,u,\nabla\varphi)v$ under the integral sign on the right hand side in (\ref{15}) is in $L^1(Q_T)$. Indeed, (H3) gives $f(\cdot,\cdot,u,\nabla\varphi)\in L^1(Q_T)$ while $v\in L^\infty\big(0,T;W^{1,q'}(\Omega)\big)$ can be identified with a function in $L^\infty(Q_T)$. To see this, we take any $r\ge q'$ and obtain 
\[
 \int\limits_\Omega\big|v(x,t)\big|^rdx\le\gamma_0^r\big\|v(t)\big\|_{W^{1,q'}}^r\operatorname{mes}\Omega\quad\text{ for a.e. } \ t\in[0,T],
\]
where $\gamma_0$ denotes the embedding constant of $W^{1,q'}(\Omega)\subset C(\overline{\Omega})$ (notice $q'> n+2$). Thus, 
\[
 \|v\|_{L^\infty(Q_T)} =\lim\limits_{r\to\infty}\left(\;\int\limits_{Q_T}\big|v(x,t)\big|^rdxdt\right)^{1/r}\le\gamma_0\|v\|_{L^\infty (0,T;W^{1,q'})}<+\infty.
\]

To make precise the meaning of (\ref{17}), let $\frac{2n}{n+2}<q<\frac{n+2}{n+1}$ (cf. our main result below). We obtain the dense and continuous embeddings $W^{1,q}(\Omega)\subset L^{nq/(n-q)}(\Omega)\subset L^2(\Omega)$. Identifying $L^2(\Omega)$ with its dual space, it follows $L^2(\Omega)\subset\big(W^{1,q'}(\Omega)\big)^*$ continuously (since $W^{1,q'}(\Omega)\subset W^{1,q}(\Omega)$ continuously). Thus,
\[
 u\in L^q\big(0,T;\big(W^{1,q'}(\Omega)\big)^*\big),\quad u'\in L^1\big(0,T;\big(W^{1,q'}(\Omega)\big)^*\big)\quad (\text{cf. (\ref{13})}).
\]
Hence, there exists $\tilde{u}\in C\big([0,T];\big(W^{1,q'}(\Omega)\big)^*\big)$ such that $\tilde{u}(t)=u(t)$ for a.e. $t\in[0,T]$. Then (\ref{17}) with initial datum (\ref{12}) is meant in the sense
\[
 \big\langle\tilde{u}(0),z\big\rangle_{W^{1,q'}}=\int\limits_\Omega u_0(x)z(x)dx\quad\forall\; z\in W^{1,q'}(\Omega).
\]

\begin{remark}\label{r2}
 Let $(\varphi,u)$ be a sufficiently regular classical solution of (\ref{1})--(\ref{5}). More specifically, let $u\in C^1(\overline{Q}_T)$. Then the function $t\mapsto u(\cdot,t)$ possesses a distributional derivative $u'\in L^2\big(0,T;L^2(\Omega)\big)$ and there holds
 \[
  \int\limits_0^T(u',v)_{L^2}dt=\int\limits_{Q_T}\frac{\partial u}{\partial t}\, v\, dxdt\quad\forall\; v\in L^2\big(0,T;L^2(\Omega)\big).
  \]

By routine arguments one obtains that $(\varphi,u)$ satisfies (\ref{14}) and (\ref{15}). Thus, (\ref{13})--(\ref{17}) represents a weak formulation of (\ref{1})--(\ref{5}).\hfill $\square$
\end{remark}

The main result of our paper is the following

\begin{main}
 Assume {\rm (H1)--(H3)}. In addition to {\rm (H1)}, suppose that
 \begin{equation}\label{18}
  \big(\sigma(u,\tau)\tau-\sigma(u,\overline{\tau})\overline{\tau}\big)(\tau-\overline{\tau})>0\quad\forall\; u\in\mathbb{R}, \ \forall\; \tau,\overline{\tau}\in\mathbb{R}_+, \; \tau\ne\overline{\tau}.
 \end{equation}
Further, let {\rm (\ref{10})} and {\rm (\ref{12})} be satisfied, and let  
\begin{equation}\label{19}
 g=\mathrm{const}>0, \quad h=\mathrm{const}. 
\end{equation}
$($cf. {\rm (\ref{11}))}. Then there exists a pair
\[
 (\varphi,u)\in L^p\big(0,T;W^{1,p}(\Omega)\big)\times\left(\bigcap\limits_{1< q<(n+2)/(n+1)}L^q\big(0,T;W^{1,q}(\Omega)\big)\right)
\]
such that {\rm (\ref{14})} and {\rm (\ref{16})} hold, and {\rm (\ref{13}), (\ref{15})} and {\rm (\ref{17})} hold for every $1<q<\frac{n+2}{n+1}$. Moreover, $u$ satisfies 
\begin{equation}\label{20}
       {\displaystyle\|u\|_{L^\infty(0,T;L^1)} +\lambda\int\limits_{Q_T}\frac{|\nabla u|^2}{(1+|u|)^{1+\lambda}}\,dxdt}
       \le c\big(1+\|u_0\|_{L^1}+\|\varphi_D\|_{L^p (0,T; W^{1,p})}^p\big);
\end{equation}
\begin{equation}
\label{21}
u\in\bigcap\limits_{1<r<(n+2)/n}L^r\big(0,T;L^r(\Omega)\big).
\end{equation}
\end{main}
\smallskip

\begin{remark}\label{r3}
 (Cf. (\ref{18})). For $a:\mathbb{R}_+\to\mathbb{R}_+$ the following two statements about strict monotonicity are equivalent
 \[
  \begin{array}{ll}
  \text{(i)}&\quad \big(a(\tau)\tau-a(\overline{\tau})\overline{\tau}\big)(\tau-\overline{\tau})>0\quad\forall\; \tau,\overline{\tau}\in\mathbb{R}_+,\;\tau\ne\overline{\tau};\\[3mm]
  \text{(ii)}&\quad\big(a\big(|\xi|\big)\xi-a\big(|\overline{\xi}|\big)\overline{\xi}\big)\cdot(\xi-\overline{\xi})>0\quad\forall\;\xi,\overline{\xi}\in\mathbb{R}^n,\;\xi\ne\overline{\xi}\;\; (n\ge 2).
  \end{array}
 \]
This can be easily verified by elementary calculations.\par
We notice that the strict monotonicity of the function $\xi\mapsto\sigma_0(u)\big(\delta+|\xi|^2\big)^{(p-2)/2}\xi$ ($\delta>0$, $1<p<+\infty$; cf. (\ref{6})) [as well as of the function $\xi\mapsto\sigma_0(u)|\xi|^{p-2}\xi$ ($2\le p<+\infty$)] follows from the inequalities
\begin{align*}
 &\big(\big(\delta+|\xi|^2\big)^{(p-2)/2}\xi-\big(\delta+|\overline{\xi}|^2\big)^{(p-2)/2}\overline{\xi}\big)\cdot(\xi-\overline{\xi})\\[1mm]
 &\ge\left\{\begin{array}{l}
             {\displaystyle\frac{p-1}{\big(\delta_0+|\xi|^2+|\overline{\xi}|^2\big)^{(2-p)/2}}|\xi-\overline{\xi}|^2\quad\forall\;0\le\delta\le\delta_0,\;\;\forall\; 1<p\le 2,}\\[5mm]
             {\displaystyle\min\Big\{\frac12,\frac1{2^{p-2}}\Big\}|\xi-\overline{\xi}|^p\quad\forall\;\delta\ge 0,\;\;\forall\; 2\le p<+\infty}
            \end{array}\right.
\end{align*}
(cf. \cite[(I) p.71, p. 74]{19}). The coefficient $\min\Big\{\frac12,\frac1{2^{p.2}}\Big\}$ is related to the inequality
\[
 \alpha^{p-2}+\beta^{p-2}\ge\min\Big\{1,\frac1{2^{p-3}}\Big\}(\alpha+\beta)^{p-2}\quad\forall\;\alpha,\beta\in\mathbb{R}_+,\;\forall\;2\le p<+\infty
\]
(note by P.-A. Ivert).
\end{remark}

\begin{remark}\label{r4}
 For $\sigma(u,\tau)=\sigma_0(u)$ (i.e., $p=2$ in (\ref{6})), $f(x,t,u,\xi)=\sigma_0(u)|\xi|^2$ (cf. (\ref{9})) and Dirichlet boundary conditions, in \cite{1} ($n=3$) and \cite{9} ($n=2$) the authors proved the existence of a weak solution $(\varphi,u)$ of (\ref{1})--(\ref{5}) such that
 \[
  (\varphi,u)\in L^2\big(0,T;W^{1,2}(\Omega)\big)\times L^2\big(0,T;W^{1,2}(\Omega)\big).
 \]
\end{remark}
\medskip

\section{Proof of the main result}

We divide the proof into three parts.

\subsection{Existence of an approximate solution}

For $\ep>0$ we define the Carath\'eodory function
\[
 f_\ep(x,t,u,\xi)=\frac{f(x,t,u,\xi)}{1+\ep f(x,t,u,\xi)},\quad (x,t,u,\xi)\in Q_T\times\mathbb{R}\times\mathbb{R}^n.
\]

Let $(u_{0,\ep})_{\ep>0}$ be a sequence of functions in $L^2(\Omega)$ such that $u_{0,\ep}\to u_0$ strongly in $L^1(\Omega)$ as $\ep\to 0$. We have
\medskip

\begin{lemma}\label{l1}
 For every $\ep>0$ there exists a pair
 \[
  (\varphi_\ep,u_\ep)\in L^p\big(0,T;W^{1,p}(\Omega)\big)\times L^2\big(0,T;W^{1,2}(\Omega)\big)
 \]
such that
\begin{equation}\label{22}
\int\limits_{Q_T}\sigma\big(u_\ep,|\nabla\varphi_\ep|\big)\nabla\varphi_\ep\cdot\nabla\zeta\, dxdt=0\quad\forall\;\zeta\in L^p\big(0,T;W_{\Gamma_D}^{1,p}(\Omega)\big);
\end{equation}
\begin{equation}\label{23}
 \varphi_\ep=\varphi_D\quad\text{a.e. on } \ \Sigma_D;
\end{equation}

\begin{equation}\label{24}
 \exists\;u'_\ep\in L^2\big(0,T;\big(W^{1,2}(\Omega)\big)^*\big);
\end{equation}
\begin{equation}\label{25}
 \left.\begin{array}{l}
        {\displaystyle\int\limits_0^T\langle u'_\ep,v\rangle_{W^{1,2}}dt+\int\limits_{Q_T}\kappa(u_\ep)\nabla u_\ep\cdot\nabla v\, dxdt+g\int\limits_0^T\int\limits_{\partial\Omega}(u_\ep-h)v\, d_xSdt}\\[6mm]
        {\displaystyle\qquad\qquad=\int\limits_{Q_T}f_\ep(x,t,u_\ep,\nabla\varphi_\ep)v\,dxdt\quad\forall\; v\in L^2\big(0,T;W^{1,2}(\Omega)\big);}
       \end{array}\right\}
\end{equation}
\begin{equation}\label{26}
 u_\ep(\cdot,0)=u_{0,\ep}\quad\text{ a.e. in } \ \Omega.
\end{equation}
\end{lemma}
\medskip

By routine arguments it is readily seen that (\ref{25}) is equivalent to 
\begin{equation}\label{27}
 \left.\begin{array}{l}
        {\displaystyle\big\langle u'_\ep(t),z\big\rangle_{W^{1,2}}+\int\limits_\Omega\kappa\big(u_\ep(x,t)\big)\nabla u_\ep(x,t)\cdot\nabla z(x)\,dx}\\[6mm]
        {\displaystyle+g\int\limits_{\partial\Omega}\big(u_\ep(x,t)-h\big)z(x)d_xS}\\[6mm]
        {\displaystyle =\int\limits_\Omega f_\ep\big(x,t,u_\ep(x,t),\nabla\varphi_\ep(x,t)\big)z(x)dx}
       \end{array}\right\}
\end{equation}
for a.e. $t\in[0,T]$ and all $z\in W^{1,2}(\Omega)$, where the set of measure zero of those $t$ for which (\ref{27}) fails, does not depend on $z$.
\bigskip

\noindent
{\it Proof of Lemma~{\rm \ref{l1}}. } \
 We prove this lemma by the aid of Schauder's Fixed Point Theorem.
 \medskip
 
 \noindent
 {\em Step $1$. Construction of a mapping}
 \[
  \mathcal{T}:\overline{\mathcal{B}}_R\to\overline{\mathcal{B}}_R,
 \]
where 
\[
 \overline{\mathcal{B}}_R=\left\{w\in L^2\big(0,T;L^2(\Omega)\big); \|w\|_{L^2(L^2)}\le R\right\}\footnotemark
\]
($R>0$ suitable choosen). For this we need the following two preliminary results $1^\circ$ and $2^\circ$.\footnotetext{In what follows, for indexes we write $L^p(L^p)$ in place of $L^p\big(0,T;L^p(\Omega)\big)$, $L^p(W_D^{1,p})$ etc.}
\smallskip

$1^\circ$ {\it Given any $u\in L^2\big(0,T;L^2(\Omega)\big)$, there exists exactly one $\varphi\in L^p\big(0,T;W^{1,p}(\Omega)\big)$ \ $(\varphi=\varphi_u)$ such that}
\begin{equation}\label{28}
 \int\limits_{Q_T}\sigma\big(u,|\nabla\varphi|\big)\nabla\varphi\cdot\nabla\zeta\, dxdt=0\quad\forall\;\zeta\in L^p\big(0,T;W_{\Gamma_D}^{1,p}(\Omega)\big);
\end{equation}
\begin{equation}\label{29}
 \varphi=\varphi_D\quad\text{ \it a.e. on } \ \Sigma_D.
\end{equation}
\smallskip

\noindent 
To prove this, we define a mapping
\[
 \mathcal{A}:L^p\big(0,T;W_{\Gamma_D}^{1,p}(\Omega)\big)\longrightarrow L^{p'}\big(0,T;\big(W_{\Gamma_D}^{1,p}(\Omega)\big)^*\big)\qquad(\mathcal{A}=\mathcal{A}_u)
\]
by
\[
 \big\langle\mathcal{A}(\psi),\zeta\big\rangle_{L^p(W_{\Gamma_D}^{1,p})}=\int\limits_{Q_T}\sigma\big(u,\big|\nabla(\psi+\varphi_D)\big|\big)\nabla(\psi+\varphi_D)\cdot\nabla\zeta\, dxdt,
\]
where $\psi,\zeta\in L^p\big(0,T;W_{\Gamma_D}^{1,p}(\Omega)\big)$. From (H1) and (\ref{10}) it follows that this mapping maps bounded sets into bounded sets. By (\ref{18}) [cf. also Remark~\ref{r4}] we have 
\[
 \big\langle\mathcal{A}(\psi)-\mathcal{A}(\overline{\psi}),\psi-\overline{\psi}\big\rangle_{L^p(W_{\Gamma_D}^{1,p})}>0
\]
for all $\psi,\overline{\psi}\in L^p(0,T;W_{\Gamma_D}^{1,p})$, $\psi\ne\overline{\psi}$. Finally, appealing once more to (H1) we obtain the coercivity of $\mathcal{A}$.\par
The theory of monotone operators yields the existence and uniqueness of an $\omega\in L^p\big(0,T;W_{\Gamma_D}^{1,p}(\Omega)\big)$ such that
\[
 \mathcal{A}(\omega)=0
\]
(see, e.g., \cite[Chap. 2.2]{20}, \cite[Chap. 26.2]{29}). Then the function $\varphi=\omega+\varphi_D$ is in $L^p\big(0,T;W^{1,p}(\Omega)\big)$ and solves (\ref{28})--(\ref{29}).
\medskip

$2^\circ$ {\it Let $u\in L^2\big(0,T;L^2(\Omega)\big)$ and let $\varphi=\varphi_u\in L^2\big(0,T;W^{1,2}(\Omega)\big)$ denote the uniquely determined solution of {\rm (\ref{28})--(\ref{29})} $($cf. $1^\circ)$. Then there exists exactly one 
\[
 \hat{u}\in L^2\big(0,T;W^{1,2}(\Omega)\big)\cap C\big([0,T];L^2(\Omega)\big)\quad (\hat{u}=\hat{u}_{\varphi_u})
\]
such that
\begin{equation}\label{30}
 \exists\; \hat{u}'\in L^2\big(0,T;\big(W^{1,2}(\Omega)\big)^*\big);
\end{equation}
\begin{equation}\label{31}
 \left.\begin{array}{l}
        {\displaystyle\int\limits_0^T\big\langle\hat{u}',v\big\rangle_{W^{1,2}}dt+\int\limits_{Q_T}\kappa(u)\nabla\hat{u}\cdot\nabla v\,dxdt+g\int\limits_0^T\int\limits_{\partial\Omega}\big(\hat{u}-h\big)v\,d_xSdt}\\[6mm]
        {\displaystyle=\int\limits_{Q_T}f_\ep(x,t,u,\nabla\varphi) v\,dxdt\quad\forall\; v\in L^2\big(0,T;W^{1,2}(\Omega)\big);}
       \end{array}\right\}
\end{equation}
\begin{equation}\label{32}
 \hat{u}(\cdot,0)=u_{0,\ep}\quad\text{ \it a.e. in } \ \Omega;
\end{equation}\vspace*{-3mm}
\begin{equation}\label{33}
 \|\hat{u}\|_{L^\infty(L^2)}+\|\hat{u}\|_{L^2(W^{1,2})}+\|\hat{u}'\|_{L^2((W^{1,2})^*)}\le c,
\end{equation}
\vspace*{-2mm}

\noindent
where the constant $c$ depends on $\kappa_0$ $\kappa_1$, $g$, $h$ $($see {\rm (H2), (\ref{19}))}, $\|u_{0,\ep}\|_{L^2}$ and $\frac1\ep$, but is independent of $u$. }
\medskip

\noindent
This result follows from the theory of linear evolution equations (see, e.g., \cite[Chap. 7.1]{12}). To see this, it suffices to notice that 
\[
 [w,z]_{W^{1,2}}=\int\limits_\Omega\nabla w\cdot\nabla z\,dx+\int\limits_{\partial\Omega} wz\,d_xS,\quad w,z\in W^{1,2}(\Omega)
\]
is a scalar product on $W^{1,2}(\Omega)$ which is equivalent to the standard scalar product on this space.\hfill $\square$
\medskip

>From (\ref{33}) we conclude that there exists a constant $R>0$ which depends the same quantities as the constant $c$ such that $\|\hat{u}\|_{L^2(L^2)}\le R$. We now define a mapping
\[
 \mathcal{T}:\overline{\mathcal{B}}_R\longrightarrow\overline{\mathcal{B}}_R
\]
by
\[
 \mathcal{T} u=\hat{u},\qquad \hat{u}\quad\text{according to $2^\circ$}.  
\]
{}\hfill$\square$
\medskip

\noindent
{\it Step $2$. Properties of $\mathcal{T}$\/} \ \ We have
\medskip

$3^\circ$ $\mathcal{T}(\overline{\mathcal{B}}_R)$ {\it is precompact;}
\medskip

$4^\circ$ $\mathcal{T}$ {\it is continuous.}
\bigskip

\noindent
{\it Proof of $3^\circ$ }  Let $(w_k)\subset\mathcal{T}(\overline{\mathcal{B}}_R)$ ($k\in\mathbb{N}$) be any sequence. Then $w_k=\mathcal{T} u_k=\hat{u}_k$, where $u_k\in\overline{\mathcal{B}}_R$. By (\ref{33}), 
\[
 \|w_k\|_{L^2(W^{1,2})}+\|w'_k\|_{L^2((W^{1,2})^*)}\le c\quad\forall\; k\in \mathbb{N}.
\]
The embedding $W^{1,2}(\Omega)\subset L^2(\Omega)$ being dense and compact, a well-known compactness theorem (see \cite[pp. 58--59]{20}) yields the existence of a subsequence of $(w_k)$ (not relabelled) such that
\[
 w_k\to w\quad\text{ strongly in }\quad L^2\big(0,T;L^2(\Omega)\big)\quad\text{as }\; k\to \infty.
\]

\noindent
{\it Proof of $4^\circ$ } Let $(u_k)\subset\overline{\mathcal{B}}_R$ ($k\in\mathbb{N}$) be a sequence such that $u_k\to u$ strongly in $L^2\big(0,T;L^2(\Omega)\big)$ as $k\to\infty$. By passing to a subsequence if necessary, we may assume
\begin{equation}\label{34}
 u_k\to u\quad\text{ a.e. in}\quad Q_T\quad\text{as } \ k\to \infty. 
\end{equation}

Let $\varphi_k,\varphi\in L^p\big(0,T;W^{1,p}(\Omega)\big)$ ($\varphi_k=\varphi_{u_k}$, $\varphi=\varphi_u$; $k\in\mathbb{N}$) be determined by $1^\circ$, i.e.,
\begin{equation}\label{35}
 \int\limits_{Q_T}\sigma\big( u_k,|\nabla\varphi_k|\big)\nabla\varphi_k\cdot\nabla\zeta\, dxdt=0\quad\forall\;\zeta\in L^p\big(0,T;W_{\Gamma_D}^{1,p}(\Omega)\big),
\end{equation}
\begin{equation}\label{36}
 \varphi_k=\varphi_D\quad\text{ a.e. on } \ \Sigma_D,
\end{equation}
and $\varphi$ satisfies (\ref{28})--(\ref{29}).\par
Next, let $\hat{u}_k,\hat{u}\in L^2\big(0,T;W^{1,2}(\Omega)\big)$ ($k\in\mathbb{N}$) be determined by $2^\circ$, i.e., $\hat{u}_k$ satisfies (\ref{30})--(\ref{33}) in place of $\hat{u}$. We show
\begin{equation}\label{37}
 \mathcal{T}u_k=\hat{u}_k\longrightarrow\hat{u}=\mathcal{T}u\quad\text{strongly in } \ L^2\big(0,T;L^2(\Omega)\big)\;\text{ as }\; k\to\infty.
\end{equation}
To this end, let us assume
\begin{equation}\label{38}
 \nabla\varphi_k\longrightarrow\nabla\varphi\quad\text{ a.e. in}\quad Q_T\;\text{ as }\; k\to\infty
\end{equation}
(the proof will be given below). We insert $v=\hat{u}_k-\hat{u}$ into the variational identities in (\ref{31}) for $\hat{u}_k$ and $\hat{u}$, respectively, and form the difference of both identities. This gives an integral relation which contains the term
\[
 \int\limits_0^T\big\langle\hat{u}_k'-\hat{u}',\hat{u}_k-\hat{u}\big\rangle_{W^{1,2}}dt=\frac12\big\|\hat{u}_k(T)-\hat{u}(T)\big\|_{L^2}^2
\]
(observe (\ref{32})) and the right hand side 
\[
 \int\limits_{Q_T}\big(f_\ep(x,t,u_k,\nabla\varphi_k)-f_\ep(x,t,u,\nabla\varphi)\big)\big(\hat{u}_k-\hat{u}\big)dxdt.
\]
By (\ref{34}) and (\ref{38}),
\[
 \lim\limits_{k\to\infty}\int\limits_{Q_T}\big(f_\ep(x,t,u_k,\nabla\varphi_k)-f_\ep(x,t,u,\nabla\varphi)\big)^2dxdt=0.
\]
The claim (\ref{37}) is proved.
\medskip

\noindent
{\it Proof of\/} (\ref{38}) \,\, By (\ref{36}), the function $\zeta=\varphi-\varphi_D$ is admissible in (\ref{35}). Combining (H1) and H\"older's inequality  we obtain $\|\varphi_k\|_{L^p(W^{1,p})}\le\mathrm{const}$ for all $k\in\mathbb{N}$. Hence, there exists a subsequence of $(\varphi_k)$ (not relabelled) such that 
\[
 \varphi_k\longrightarrow\chi\quad \text{ weakly in}\quad L^p\big(0,T;W^{1,p}(\Omega)\big)\quad\text{as }\; k\to\infty.
\]
It follows $\chi=\varphi_D$ a.e. on $\Sigma_D$. Observing (\ref{34}), the passage to the limit $k\to\infty$ in (\ref{35}) is easily carried out by the monotonicity trick (see, e.g., \cite[p. 172]{20}, \cite[p.~474]{29}) to obtain
\[
 \int\limits_{Q_T}\sigma\big(u,|\nabla\chi|\big)\nabla\chi\cdot\nabla\zeta\, dxdt=0\quad\forall\; \zeta\in L^p\big(0,T;W_{\Gamma_D}^{1,p}(\Omega)\big).
\]
Thus, $\chi$ satisfies (\ref{28})--(\ref{29}) in place of $\varphi$. By the strict monotonicity of $\xi\mapsto\sigma\big(u,|\xi|\big)\xi$ (cf. (\ref{18}) resp. Remark~\ref{r3}) we obtain $\chi=\varphi$, and the whole sequence $(\varphi_k)$ converges weakly in $L^p\big(0,T;W^{1,p}(\Omega)\big)$ to $\varphi$. Therefore, 
\[
 \lim\limits_{k\to\infty}\int\limits_{Q_T}\big[\sigma\big(u_k,|\nabla\varphi_k|\big)\nabla\varphi_k-\sigma\big(u_k,|\nabla\varphi|\big)\nabla\varphi\big]\cdot\nabla(\varphi_k-\varphi)dxdt=0.
\]
Finally, for a.e. $(x,t)\in Q_T$, define
\begin{align*}
 E_k(x,t)&=\big[\sigma\big(u_k(x,t),\big|\nabla\varphi_k(x,t)\big|\big)\nabla\varphi_k(x,t)\\
 &\quad -\sigma\big(u_k(x,t),\big|\nabla\varphi(x,t)\big|\big)\nabla\varphi(x,t)\big]\cdot\nabla\big(\varphi_k(x,t)-\varphi(x,t)\big).
\end{align*}
We obtain
\[
E_k(x,t)\longrightarrow 0\quad\text{ for a.e.}\quad (x,t)\in Q_T\;\text{ as }\; k\to\infty.
\]
A well-known argument due to Leray-Lions \cite{17} now gives (\ref{38}) (cf. also, \cite[pp.~184--185]{20}, \cite{2}, \cite{3}, \cite{23}, \cite{24}).
\medskip

\noindent
{\it Step $3$. \ Existence of a fixed point of $\mathcal{T}$\/} \,\, The Schauder Fixed Point Theorem yields the existence of an element $u_\ep\in\overline{\mathcal{B}}_R$ such that $\mathcal{T} u_\ep=u_\ep$. We then determine $\varphi_\ep=\varphi_{u_\ep}$ according to $1^\circ$ above. The pair $(\varphi_\ep,u_\ep)$ satisfies (\ref{22})--(\ref{26}).

\subsection{A-priori estimates}

We have

\begin{lemma}\label{l2}
Let be $(\varphi_\ep,u_\ep)$ as in Lemma~{\rm \ref{l1}}. 
\begin{equation}\label{39}
 \|\varphi_\ep\|_{L^p(W^{1,p})}\le c\big(1+\|\varphi_D\|_{L^p(W^{1,p})}\big)\footnotemark;
\end{equation}
\footnotetext{In what follows, by $c$ we denote constants which may change their numerical value from line to line but do not depend on $\ep$.}
\begin{equation}\label{40}
 \int\limits_{Q_T} f_\ep(x,t,u_\ep,\nabla\varphi_\ep)dxdt\le c\big(1+\|\varphi_D\|_{L^p(W^{1,p})}^p\big);
\end{equation}
\begin{equation}\label{41}
 \left.\begin{array}{l}
        {\displaystyle\|u_\ep\|_{L^\infty(L^1)}+\lambda\int\limits_{Q_T}\frac{|\nabla u_\ep|^2}{\big(1+|u_\ep|\big)^{1+\lambda}}dxdt}\\[6mm]
        {\displaystyle\le c\big(1+\|u_{0,\ep}\|_{L^1}+\|\varphi_D\|_{L^p(W^{1,p})}^p\big),\quad 0<\lambda<1;}
       \end{array}\right\}
\end{equation}
\begin{align}\label{42}
 \|u_\ep\|_{L^q(W^{1,q})}\le& \;c(q)\quad\forall\; 1<q<{\displaystyle\frac{n+2}{n+1}};\\
 \label{43}
 \| u_\ep\|_{L^r(L^r)}\le& \;c(r)\quad\forall\; 1<r<\frac{n+2}n;\\
\label{44}
 \|u_\ep'\|_{L^1((W^{1,q'})^*)}\le &\;c(q)\quad \forall\; 1<q<\frac{n+2}{n+1},
\end{align}
\vspace*{-2mm}

\noindent
where $c(q)\to+\infty$ as $q\to\frac{n+2}{n+1}$, and $c(r)\to+\infty$ as $r\to\frac{n+2}n$.
\end{lemma}
\medskip

\noindent
{\it Proof of Lemma~{\rm \ref{l2}}. }
 By (\ref{23}) the function $\zeta=\varphi_\ep-\varphi_D$ is admisible in (\ref{22}). The estimate (\ref{39}) is then easily obtained by (H1) and H\"older's inequality. From (H1) and (H3) it follows that 
 \begin{equation}\label{45} 
 f(x,t,u,\xi)\le c_5\big(1+\sigma\big(u,|\xi|\big)|\xi|^2\big)\quad (c_5=\mathrm{const}>0)
 \end{equation}
 \vspace*{-3mm}
 
 \noindent
for all $(x,t,u,\xi)\in Q_T\times\mathbb{R}\times\mathbb{R}^n$. With the help of this inequality the estimate in (\ref{40}) is easily deduced from (\ref{39}).\par
To prove (\ref{41}), for $s\in\mathbb{R}$ and $0<\lambda<1$ we define the functions
\begin{align*}
 \Phi(s)&=\Phi_\lambda(s)=\left(1-\frac1{\big(1+|s|\big)^\lambda}\right)\operatorname{sign} s\:\footnotemark,\\
 \Psi(s)&=\Psi_\lambda(s)=|s|+\frac1{1-\lambda}\big(1-\big(1+|s|\big)^{1-\lambda}\big).                                                                                \end{align*}
\footnotetext{$\operatorname{sign}(0)=0$.}
We obtain 
\[
\Phi'(s)=\frac\lambda{\big(1+|s|\big)^{1+\lambda}},\quad\Psi'(s)=\Phi(s),\quad \frac{|s|}2-\frac{2^{(1-\lambda)/2}}{1-\lambda}\le\Psi(s)\le|s|,           
\]
and
\[
 \nabla\Phi(u_\ep)=\Phi'(u_\ep)\nabla u_\ep=\lambda\frac{\nabla u_\ep}{\big(1+|u_\ep|\big)^{1+\lambda}}\quad\text{ for a.e. }\; (x,t)\in Q_T,
\]
\[
 \int\limits_0^t\big\langle u'_\ep,\Phi(u_\ep)\big\rangle_{W^{1,2}}ds=\int\limits_\Omega\Psi\big(u_\ep(x,t)\big)dx-\int\limits_\Omega\Psi\big(u_{0,\ep}(x)\big)dt\quad\forall\; t\in [0,T]
\]
(cf. \cite{22}, \cite{23}, \cite{24}). We insert $z=\Phi\big(u_\ep(\cdot,t)\big)$ into (\ref{27}), integrate over the interval $[0,t]$ and make use of (\ref{40}). By elementary calculations we obtain (\ref{41}).\par
The proof of (\ref{42}) is now easily done by well-known arguments (see, e.g., \cite{2}, \cite{3}, \cite{23}, \cite{24}). Indeed, let $1<q<n$. A simple application of H\"older's inequality yields
\begin{equation}\label{46}
 \|z\|_{L^{q(n+1)/n}}\le\|z\|_{L^1}^{1/(n+1)}\|z\|_{L^{nq/(n-q)}}^{n/(n+1)}\quad\forall\; z\in L^{nq/(n-q)}(\Omega).
\end{equation}
Next, given $1<q<\frac{n+2}{n+1}$ we set $\lambda=\frac 1n\big(n+2-q(n+1)\big)$. Using the integral estimate in (\ref{41}) one finds
\begin{align*}
 \int\limits_{Q_T}|\nabla u_\ep|^qdxdt&=\int\limits_{Q_T}\frac{|\nabla u_\ep|^q}{\big(1+|u_\ep|\big)^{q(1+\lambda)/2}}\cdot\big(1+|u_\ep|\big)^{q(1+\lambda)/2}dxdt\\
 &\le\frac c{\lambda^{q/2}}\left(\;\int\limits_{Q_T}\big(1+|u_\ep|\big)^{q(n+1)/n}dxdt\right)^{(2-q)/2}.
\end{align*}
To estimate the integral on the right hand side, we take $z=u_\ep(\cdot,t)$ in (\ref{46}) and use then the bound on $\|u_\ep\|_{L^\infty(L^1)}$ in (\ref{41}) and the Sobolev embedding theorem $W^{1,q}(\Omega)\subset L^{nq/(n-q)}(\Omega)$. We obtain 
\[
 \int\limits_{Q_T}\big(1+|u_\ep|\big)^{q(n+1)/n}dxdt\le c\left(1+\int\limits_{Q_T}|\nabla u_\ep|^qdxdt\right).
\]
Whence (\ref{42}). Using once more (\ref{46}), we get (\ref{43}).
\smallskip

We finally prove (\ref{44}). From (\ref{27}) it follows for a.e. $t\in [0,T]$ and all $z\in W^{1,q'}(\Omega)$ that
\[
 \big|\big\langle u'_\ep(t),z\big\rangle_{W^{1,q'}}\big|\le c\big(1+\big\|u_\ep(t)\big\|_{W^{1,q}}+\big\|\varphi_\ep(t)\big\|_{W^{1,p}}^p\big)\|z\|_{W^{1,q'}},
\]\vspace*{-3mm}

\noindent
where the constant $c$ does not depend on $\ep$. By (\ref{39}) and (\ref{40}), the function in parantheses is uniformly bounded independently of $\ep$. The estimate (\ref{44}) is now easily seen.

\subsection{Passage to the limit $\ep\to 0$}

We begin by proving the existence of convergent subsequences of 
$(\varphi_\ep,u_\ep)$. Then we complete the proof of our main result by showing that the limit functions of these subsequences yield a weak solution of (\ref{1})--(\ref{5}).

\begin{lemma}\label{l3}
Let be $(\varphi_\ep,u_\ep)$ as in Lemma~{\rm \ref{l1}}. Then there exists a subsequence $($not relabelled$)$ such that
\begin{equation}\label{47}
\varphi_\ep\longrightarrow\varphi\quad\text{ weakly in}\quad L^p\big(0,T;W^{1,p}(\Omega)\big); 
 \end{equation}
\begin{equation}\label{48}
 \left.\begin{array}{l}
        {\displaystyle u_\ep\longrightarrow u\quad\text{ weakly in}\quad L^q\big(0,T;W^{1,q}(\Omega)\big)\quad \Big(1<q<\frac{n+2}{n+1}\Big)}\\[4mm]
        {\displaystyle \text{and weakly in}\quad L^r\big(0,T;L^r(\Omega)\big)\quad \Big(1<r<\frac{n+2}n\Big);}
       \end{array}\right\}
\end{equation}
\begin{equation}\label{49} 
u_\ep\longrightarrow u\quad\text{ a.e. in}\quad Q_T;
\end{equation}
\begin{equation}\label{50}
 \nabla\varphi_\ep\longrightarrow\nabla\varphi\quad\text{ a.e. in}\quad Q_T;
\end{equation}
\begin{equation}\label{51}
\int\limits_{Q_T}\left|\sigma\big( u_\ep,|\nabla\varphi_\ep|\big)|\nabla\varphi_\ep|^2-\sigma\big(u,|\nabla\varphi|\big)|\nabla\varphi|^2\right|dxdt\longrightarrow 0 
\end{equation}
as $\ep\to 0$.
\end{lemma}

\noindent
{\it Proof of Lemma~{\rm \ref{l3}}. } The existence of subsequences of $(\varphi_\ep, u_\ep)$ satisfying (\ref{47}), (\ref{48}) follows from the reflexivity of the respective spaces.\par
We prove (\ref{49}). To this end, take $q$ such that $\frac2{n+2}<q<\frac{n+2}{n+1}$. Then $W^{1,q}(\Omega)\subset L^2(\Omega)$ compactly (recall $n=2$ or $n=3$). We identify $L^2(\Omega)$ with its dual space and obtain the continuous embedding $L^2(\Omega)\subset\big(W^{1,q'}(\Omega)\big)^*$ (cf. Section~2). Observing the bounds on $u_\ep$ and $u'_\ep$ in (\ref{42}) and (\ref{44}), respectively, from the compactness result in \cite[Prop. 1]{6} resp. \cite[Cor. 4]{26} we obtain $u_\ep\to u$ strongly in $L^q\big(0,T; L^2(\Omega)\big)$, and thus
\[
 u_\ep\longrightarrow u\quad\text{ a.e. in}\quad Q_T\quad\text{as }\;\ep\to 0
\]
(again by passing to a subsequence if necessary). With the help of this convergence of $(u_\ep)$ we find (\ref{50}) by the same arguments as for the proof of (\ref{38}).
\smallskip

It remains to prove (\ref{51}). From (H1) and (\ref{39}) it follows that the sequence $\big(\sigma\big(u_\ep,|\nabla\varphi_\ep|\big)\nabla\varphi_\ep\big)$ is bounded in $\big[L^{p'}(Q_T)\big]^n$ for all $\ep>0$. We therefore may assume that
\[
 \sigma\big(u_\ep,|\nabla\varphi_\ep|\big)\nabla\varphi_\ep\longrightarrow F\quad\text{ weakly in}\quad \big[L^{p'}(\Omega)\big]^n \quad\text{as }\;\ep\to 0.
\]
Then (\ref{22}) implies
\[
 \lim\limits_{\ep\to 0}\int\limits_{Q_T}\sigma\big(u_\ep,|\nabla\varphi_\ep|\big)|\nabla\varphi_\ep|^2dxdt=\int\limits_{Q_T} F\cdot\nabla\varphi\, dxdt.
\]
On the other hand, for all $G\in\big[L^p(\Omega)\big]^n$ and a.e. $(x,t)\in Q_T$,
\[
 \big(\sigma\big(u_\ep,|G|\big)G-\sigma\big(u_\ep,|\nabla\varphi_\ep|\big)\nabla\varphi_\ep\big) \cdot(G-\nabla\varphi_\ep)\ge 0.
\]
Integrating this inequality over $Q_T$, letting $\ep\to 0$ and using the monotonicity trick (within the context of the dual pairing $\big( L^p(Q_T),L^{p'}(Q_T)\big)$) we get
\[
 \int\limits_{Q_T}\sigma\big(u,|\nabla\varphi|\big)|\nabla\varphi|^2dxdt=\int\limits_{Q_T} F\cdot\nabla \varphi\, dxdt.
\]
Thus,
\begin{equation}\label{52}
 \lim\limits_{\ep\to 0}\int\limits_{Q_T}\sigma\big( u_\ep,|\nabla\varphi_\ep|\big)|\nabla\varphi_\ep|^2dxdt=\int\limits_{Q_T}\sigma\big(u,|\nabla\varphi|\big)|\nabla\varphi|^2 dxdt.
\end{equation}
In addition, by (\ref{49}) and (\ref{50}),
\begin{equation}\label{53}
 \sigma\big( u_\ep,|\nabla\varphi_\ep|\big)|\nabla\varphi_\ep|^2\longrightarrow\sigma\big(u,|\nabla\varphi|\big)|\nabla\varphi|^2\quad\text{a.e. in}\quad Q_T
\end{equation}
as $\ep\to 0$. Then (\ref{51}) follows from (\ref{52}) and (\ref{53}) by the aid of Lebesgue's Dominated Convergence Theorem.\hfill$\square$
\medskip

\noindent
{\it Completion of the proof of the main result. } We prove that the pair
\[
 (\varphi,u)\in L^p\big(0,T;W^{1,p}(\Omega)\big)\times\left(\bigcap\limits_{1<q<(n+2)/(n+1)}L^q\big(0,T;W^{1,q}(\Omega)\right)
\]
obtained by Lemma~\ref{l3}, fulfills all conditions stated in our main theorem. \par
The passage to the limit $\ep\to 0$ in (\ref{22}), (\ref{23}) gives (\ref{14}), (\ref{16}), respectively. The estimates in (\ref{20}) as well as the integrability property (\ref{21}) are easily derived from (\ref{41}) and (\ref{43}), respectively, and using (\ref{49}).\par
We prove that $(\varphi,u)$ satisfies (\ref{13}), (\ref{15}) and (\ref{17}). We take $z\in W^{1,q'}(\Omega)$ ($1<q<\frac{n+2}{n+1}$) in (\ref{27}), multiply each term by $\alpha\in C^1\big([0,T]\big)$, $\alpha(T)=0$, and integrate over the interval $[0,T]$. It follows
\begin{align}
 &-\int\limits_{Q_T}u_\ep z\alpha'dxdt+\int\limits_{Q_T}\kappa(u_\ep)\nabla u_\ep\cdot\nabla z\alpha\, dxdt+g\int\limits_0^T\int\limits_{\partial\Omega}(u_\ep-h)z\alpha\, d_xSdt\nonumber\\
 \label{54}
 &=\int\limits_\Omega u_{0,\ep}(x)z(x)dx\alpha(0)+\int\limits_{Q_T} f_\ep(x,t,u_\ep,\nabla\varphi_\ep)z\alpha\, dxdt.
\end{align}
The passage to the limit $\ep\to 0$ for the second integral on the right hand side of (\ref{54}) is easily done as follows. More generally, for all $w\in L^\infty(Q_T)$ we have
\begin{equation}\label{55}
 \lim\limits_{\ep\to 0}\int\limits_{Q_T} f_\ep(x,t,u_\ep,\nabla\varphi_\ep)w\,dxdt=\int\limits_{Q_T} f(x,t,u,\nabla\varphi)w\,dxdt.
\end{equation}
Indeed, we use once more (\ref{45}) to obtain
\begin{equation}\label{56}
f(x,t,u_\ep,\nabla\varphi_\ep)|w|\le c_5\|w\|_{L^\infty}\left(1+\sigma\big(u_\ep,|\nabla\varphi_\ep|\big)|\nabla\varphi_\ep|^2\right) 
\end{equation}
for a.e. $(x,t)\in Q_T$ and all $\ep>0$. Integrating this inequality over $Q_T$ and using (\ref{39}) we find
\[
\ep f(x,t,u_\ep,\nabla\varphi_\ep)\longrightarrow 0\quad\text{ for a.e.}\quad (x,t)\in Q_T\quad\text{as }\;\ep\to 0.
\]
It follows
\[
 \lim\limits_{\ep\to 0} f_\ep(x,t,u_\ep,\nabla\varphi_\ep)w=f(x,t,u,\nabla\varphi)w\quad\text{ for a.e.}\quad (x,t)\in Q_T.
\]
In addition, since $f_\ep\le f$, (\ref{56}) implies that the function on the right hand side is an integrable bound for $f_\ep |w|$ a.e. in $Q_T$. Moreover, by (\ref{51}), these functions converge in $L^1(Q_T)$ when $\ep\to 0$. The claim (\ref{55}) thus follows from Lebesgue's Dominated Convergence Theorem. \par
The passage to the limit $\ep\to 0$ in (\ref{54}) now gives
\begin{align}
 &-\int\limits_{Q_T}u z\alpha'dxdt+\int\limits_{Q_T}\kappa(u)\nabla u\cdot\nabla z\alpha\, dxdt+ g\int\limits_0^T\int\limits_{\partial\Omega}(u-h)z\alpha\, d_xSdt\nonumber\\
 \label{57} 
 &=\int\limits_\Omega u_0(x)z(x)dx\alpha(0)+\int\limits_{Q_T} f(x,t,u,\nabla\varphi)z\alpha\, dxdt
\end{align}
(recall $u_{0,\ep}\to u_0$ strongly in $L^1(\Omega)$ as $\ep\to 0$). We prove the existence of a distributional derivative of $u$. For this, we define a mapping $F:[0,T]\to\big(W^{1,q'}(\Omega)\big)^*$ by
\begin{align*}
 \big\langle F(t),z\big\rangle_{W^{1,q'}}=&-\int\limits_\Omega\kappa\big(u(x,t)\big)\nabla u(x,t)\cdot\nabla z(x)dx-g\int\limits_{\partial\Omega}\big(u(x,t)-h\big)z(x)d_xS\\
 &+\int\limits_\Omega f\big(x,t,u(x,t),\nabla\varphi(x,t)\big)z(x)dx,\quad z\in W^{1,q'}(\Omega).
\end{align*}
It follows 
\[
 \big\|F(t)\big\|_{(W^{1,q'})^*}\le c\left(1+\big\|u(t)\big\|_{W^{1,q}}+\big\|\varphi(t)\big\|_{W^{1,p}}^p\right)
\]
for a.e. $t\in[0,T]$ (cf. the proof of (\ref{44})). The function $t\mapsto\big\langle F(t),z\big\rangle_{W^{1,q'}}$ being measurable for all $z\in W^{1,q'}(\Omega)$, $F$ is strongly measurable on $[0,T]$ by virtue of Pettis' Theorem. Thus, 
\[
 F\in L^1\big(0,T;\big(W^{1,q'}(\Omega)\big)^*\big).
\]
We take $\alpha\in C_c^\infty\big(\,]\,0,T\,[\,\big)$ in (\ref{57}) and rewrite this variational identity in the form 
\[
 \left\langle-\int\limits_0^T u(t)\alpha'(t)dt,z\right\rangle_{W^{1,q'}}=\left\langle\int\limits_0^T F(t)\alpha(t)dt,z\right\rangle_{W^{1,q'}}.
\]
This implies the existence of the distributional derivative $u'\in L^1\big(0,T;\big(W^{1,q'}(\Omega)\big)^*\big)$ and 
\[
 \int\limits_0^T u'(t)\alpha(t)dt\mathop{=}\limits^{\mathrm{in}\:(W^{1,q'})^*}\int\limits_0^T F(t)\alpha(t)dt\quad \forall\; \alpha\in C_c^\infty\big(\,]\,0,T\,[\,\big).
\]
It follows that there exists $\tilde{u}\in C\big([0,T];\big(W^{1,q'}(\Omega)\big)^*\big)$ such that $\tilde{u}(t)=u(t)$ for a.e. $t\in[0,T]$ and
\[
 \left\langle\int\limits_0^Tu'(t)\alpha(t)dt,z\right\rangle_{W^{1,q'}}=-\big\langle \tilde{u}(0)\alpha(0),z\big\rangle_{W^{1,q'}}-\int\limits_{Q_T}uz\alpha'dxdt
\]
for all $\alpha\in C^1\big([0,T]\big)$, $\alpha(T)=0$. Now (\ref{15}) and (\ref{17}) are easily obtained by standard arguments.
\medskip

\subsection*{Acknowledgement} \ The author is indebted to A. Fischer (Institute of Applied Photophysics (IAPP), Technical University Dresden) and J. Griepentrog and J. Wolf (Department of Mathematics, Humboldt University Berlin) for useful discussions when prepairing this paper.


\medskip
Received xxxx 20xx; revised xxxx 20xx.
\medskip

\end{document}